\documentclass[]{imsart}
\usepackage{amsfonts}
\usepackage[dvips]{graphics}
\usepackage{amsmath}
\usepackage{amssymb}
\usepackage{latexsym}
\RequirePackage{natbib} 
\usepackage{color}

\newtheorem{theorem}{Theorem}[section]
\newtheorem{lemma}[theorem]{Lemma}
\newtheorem{corollary}[theorem]{Corollary}
\newtheorem{proposition}[theorem]{Proposition}
\newtheorem{example}[theorem]{Example}
\newtheorem{remark}[theorem]{Remark}

\newtheorem{definition}[theorem]{Definition}
\def\bit{\begin{itemize}}
\def\eit{\end{itemize}}
\reversemarginpar   %%puts label in LEFT margin
\def\bc{\begin{center}}
\def\ec{\end{center}}
\def\bthm{\begin{theorem}}
\def\ethm{\end{theorem}}
\def\bcor{\begin{corollary}}
\def\ecor{\end{corollary}}
\def\bprop{\begin{proposition}}
\def\eprop{\end{proposition}}
\def\blem{\begin{lemma}}
\def\elem{\end{lemma}}

\def\brem{\begin{remark}}
\def\erem{\end{remark}}

\def\bdes{\begin{description}}
\def\edes{\end{description}}

\def\beq{\begin{equation}}
\def\eeq{\end{equation}}
\def\ben{\begin{enumerate}}
\def\een{\end{enumerate}}
\def\beqar{\begin{eqnarray}}
\def\eeqar{\end{eqnarray}}
\def\beqarr{\begin{eqnarray*}}
\def\eeqarr{\end{eqnarray*}}

\def\R#1{{\mathbb R}^{#1}}   %bold superscripted R-variable
\def\rar{\rightarrow}

\newcommand{\F}{\mathcal{F}}

\def\d#1dt{\frac{d#1}{dt}}    %%variable ODE left side

\def\R{\mathbb R}

\def\N{\mathbb N}

\def\hop{{\noindent}}

\newcommand{\vs}[1]{\vspace{#1}}

\newcommand{\begitem}{\begin{itemize}}
\newcommand{\finit}{\end{itemize}}

\newcommand{\zdeux}{\vs{0.2cm}} % < \medskip = 0.25cm
\newcommand{\ztrois}{\vs{0.3cm}}
\def\R{\mathbb R}
\def\N{\mathbb N}

\def\hop{{\noindent}}
\def\Argmax{\mathsf{Argmax}}
\def\rar{\rightarrow}

\begin{document}

\thispagestyle{empty}

\begin{frontmatter}

\title{Consistency of vanishingly smooth fictitious play}
\runtitle{Consistency of VSFP}

\begin{aug}
  \author{\fnms{Mathieu}  \snm{Faure}\corref{}\ead[label=e1]{mathieu.faure@univ-amu.fr}}
  \and
  \author{\fnms{Michel} \snm{Bena\"{i}m}\ead[label=e2]{sschreiber@ucdavis.edu}}

  %\thankstext{t2}{Footnote to the first author with the `thankstext' %command.}

  %\runauthor{F. Author et al.}

  \affiliation{Aix-Marseille University (Aix-Marseille School of Economics), CNRS \& EHESS \\ and \\
  Institut de Math\'ematiques, Universit\'e de Neuch\^atel,}

  \address{GREQAM, centre de la vieille charit\'e, \\
  2 rue de la vieille charit\'e \\
  13236 Marseille Cedex 02.\\
            }
          
  \address{ Rue Emile-Argand 11. Neuch\^atel. Switzerland, }
 
\end{aug}

\begin{abstract}
We discuss consistency of \emph{Vanishingly Smooth Fictitious Play}, a strategy in the context of game theory, which can be regarded as a \emph{smooth fictitious play procedure}, where the smoothing parameter is time-dependent and asymptotically vanishes. This answers a question initially raised by Drew Fudenberg and Satoru Takahashi.
\end{abstract}



\end{frontmatter}

\section{Introduction and background}

A recurring question in the theory of repeated games is to define properly a notion of \emph{good strategy} for a player facing an unknown environment. Consequently, in this paper, we are not concerned with the formalisation of  strategic interactions between rational players, but rather between a  \emph{decision maker} and \emph{nature}. Not much is known about the latter, no assumption is made on its payoff function, its thinking process or its rationality. We take the point of view of the former, whose objective is to maximize his/her average payoff in the long run. A naive approach in this direction is to assume that the game is zero-sum and to look for optimal strategies. However, the fact that his/her opponent might not try to maximize his/her payoff could lead to bad outcomes. A possible definition of \emph{good strategy} for the decision maker has been proposed by Hannan (see \cite{Han57}). It is closely related to the concept of \emph{regret}. After $n$ stages, the regret of the decision maker is the difference between the payoff that he could have obtained if he knew in advance the empirical moves of nature and the average payoff he actually got. A good strategy for the decision maker may then be defined as a strategy which ensures that, regardless of the behaviour of nature, the regret asymptotically goes to zero. Such a strategy  is called \emph{consistent}. Consistent strategies are known to exist for a long time and can be constructed, for instance, using so-called \emph{block-annealing} procedures (see e.g. \cite{Bla54}, \cite{FosVoh93}, \cite{FosVoh98} and \cite{HarMas01}). For a complete bibliography on the topic, see the last quoted paper. Also, for a recent comprehensive overview about consistency in games, see \cite{Per10} (in french).
However \emph{fictitious play} strategies are known to be non-consistent (see \cite{FudLev98}) while \emph{smooth fictitious play} strategies have been shown to be "almost" consistent by Fudenberg and Levine \cite{FudLev95} (see section \ref{ss:consistency} for a rigorous exposition).
In this paper, we consider a time-varying smooth fictitious play with a smoothing parameter decreasing to zero, that we call {\em vanishingly smooth fictitious play} (VSFP). VSFP strategies initially behave like smooth fictitious play and asymptotically  like fictitious play. The main objective of this work is to answer the following question raised to us by Drew Fudenberg and Satoru Takahashi: "are VSFP strategies consistent?"
%%%%%%%%%%%%%%%%%%%%%%%
\subsection{Notation}
\label{ss:notation}
We consider a two-player finite game in normal form. $I$ and $L$ are the (finite) set of moves of respectively player 1 (the decision maker) and player 2 (the nature).
The map $\pi: I \times L \rightarrow \mathbb{R}$ denotes the payoff function of player $1$. The sets of mixed strategies available to players are denoted $X = \Delta(I)$ and $Y = \Delta(L)$, where
\[\Delta(I) := \left\{ x \in \mathbb{R}_+^{I} \, \mid \, \; \sum_{i \in I} x_i =1 \right\},\]
and analogously  for $\Delta(L)$. As usual $\pi$ is extended to $X \times Y$ by multilinearity:
\[\forall x \in X, \, y \in Y, \, \pi(x,y) = \sum_{i \in I} \sum_{l \in L} \pi(i,l) x_i y_l.\]
\hop In the following, $(i_1,...,i_n,...)$ (respectively $(l_1,...,l_n,...)$) will denote the sequence of actions picked by player $1$ (resp.  his/her  opponent).
Let $(\Omega,\mathcal{F},\mathbb{P})$ be a probability space, endowed with a filtration $(\mathcal{F}_n)_n$. Formally, a \emph{strategy} for player $1$ is a choice of an adapted process $(i_n)_n$ on $(\Omega, (\mathcal{F}_n)_n,\mathbb{P})$, 
%i.e. a choice  of the probability distributions  $\mathbb{P}\left(i_{n+1} = \cdot \mid \mathcal{F}_n \right), \; n \in \mathbb{N}$.
%In the whole paper, we assume that  player $1$ observes his/her payoff function as well as the actions of his/her opponent. Hence, if we assume that the agents choose their next actions  according only to the past actions then a strategy for player $1$ can simply be seen as a map from $\cup_n (I \times L)^n$ to $\Delta(I)$,
%which to a given finite history $h_n  = (i_1,l_1,...,i_n,l_n)$ associates a mixed action $\sigma(h_n)$.
%i.e. the choice of a sequence of random variables $(\mathbb{P}\left(i_{n+1} = \cdot \mid \mathcal{F}_n \right))_n$.
Throughout the paper, we assume that the  agents play independently: specifically, for $(i,l) \in I \times L$, we have
\[\mathbb{P} \left(i_{n+1} = i , \; l_{n+1}=l \mid \mathcal{F}_n \right) = \mathbb{P} \left(i_{n+1} = i \mid \mathcal{F}_n \right)  \mathbb{P} \left( l_{n+1}=l \mid \mathcal{F}_n \right).\]
Finally, we call $$\overline{x}_n = \frac{1}{n} \sum_{k=1}^n \delta_{i_k}$$ the average moves of player $1$ at time $n$, $\overline{y}_n$ the  average moves of player $2$ and $$\overline{\pi}_n = \frac{1}{n} \sum_{k=1}^n \pi(i_k,l_k)$$ the average payoff to player $1.$

%%%%%%%%%%%%%%%%%%%%%%%%%%%%%%%
%%%%%%%%%%%%%%%%%%%%%%%%%%%%%%%%%%%%%%%%%%%%%%%%%%%%%%%%%%%%%%%%%%%%%%
\subsection{Consistency, definition and comments} \label{ss:consistency}
%%%%%%%%%%%%%%%%%%%%%%%%%%%%%%%%%%%%%%
We now introduce $\Pi: Y \rightarrow \mathbb{R}$, defined by
\[\Pi(y) := \max_{i \in I} \pi(i,y).\]

\hop A strategy of the decision maker is \emph{consistent} if, against any strategy of nature, the average payoff obtained by player 1 is at least as much as if the sequence of empirical moves of nature was known in advance, and decision maker had played a best response against it. More precisely, let us define the {\itshape average regret evaluation} along a sequence of moves  $h_n = (i_1,l_1,...i_n,l_n)$:
\[ e_n :=  \max_{i \in I} \pi \left(i,\frac{1}{n} \sum_{m=1}^n l_m \right) - \frac{1}{n} \sum_{m=1}^n \pi(i_m,l_m) = \Pi(\overline{y}_n) - \overline{\pi}_n.\]

\begin{definition} A strategy for player 1 is said to be consistent if, for any strategy of nature,
\[\limsup_n e_n \leq 0, \; \; \mathbb{P}-\mbox{almost surely}. \]
It is $\eta$-consistent if
\[\limsup_n e_n \leq \eta, \; \; \mathbb{P}-\mbox{almost surely}. \]
\end{definition}

Given $y \in Y$, we call $br(y)$ the set of best responses of player $1$ to $y$, namely, $$br(y) = \Argmax_{x \in X} \pi(x,y).$$ The discrete-time \emph{fictitious play} (FP) process has been introduced in \cite{Bro51}. We say that player $1$ uses a FP strategy, with prior $\overline{y}_0$ if, for $n \geq 1$, $$\mathbb{P}(i_{n+1} = \cdot \mid \mathcal{F}_n) \in br(\gamma_n),$$ where $\gamma_n = \frac{1}{n+1} \overline{y}_0 + \frac{n}{n+1} \overline{y}_n$.
It is well known that this strategy is not consistent. A simple example is given by the following (see e.g. \cite{FudLev98}).
\begin{example} \label{ex:1}
{\rm Assume that the game is matching pennies, i.e. the payoff matrix of player $1$ is given by

$$\bordermatrix{&H&T\cr
                H& 1 & 0 \cr
                T& 0 & 1}
$$
and the  prior is $\overline{y}_0 = (1/3,2/3)$. If player two acts accordingly to the deterministic rule \emph{heads} (H) on odd stages and \emph{tails} (T) on even stages, then player $1$ and $2$ always play the opposite and the average regret satisfies  $\lim_{n \rar \infty} e_n = 1/2$.}
\end{example}

However, $\eta$-consistency can be achieved by small modifications of fictitious play, which are usually called \emph{stochastic fictitious play strategies}. Originally, stochastic fictitious play was introduced by Fudenberg and Kreps in \cite{FudKre93} and the concept behind this is that players use fictitious play in a game where payoff functions are perturbed by some random variables in the spirit of  Harsanyi \cite{Har73}. On the subject, see also  \cite{FudLev95}, \cite{FudLev98} or \cite{BenHir99}.  In this paper, we adopt another point of view and assume that player $1$ chooses to randomize his/her moves by adding a small perturbation function to his/her initial payoff map $\pi$.

The class of perturbation functions usually considered (in \cite{FudLev98} or \cite{HofSan02} for instance) is the following: Consider the  maps $\rho: Int(X) \rightarrow \R$ such that:
\begin{itemize}
\item[$(A1)$] the second derivative of $\rho$ in $x$, $D^2\rho(x)$ is positive definite on the tangent space of $X$,
\[TX := \left\{h \in \R^I : \sum_i h_i = 0\right\}.\]
\item[$(A2)$] the first derivative of $\rho$, defined on $Int(X)$, verifies
\[\lim_{x \rightarrow \partial X} \|\nabla \rho(x)\| = + \infty.\]
\end{itemize}

\hop We introduce the perturbed payoff function $\tilde{\pi}$ defined, for $x \in X$,  $y \in Y$ and $\beta>0$ by
\[\tilde{\pi}(x,y,\beta) = \pi(x,y) + \frac{1}{\beta}\rho(x).\]
Under $(A1)$ and $(A2)$, the function $\tilde{\pi}$ enjoys the following property:
\begitem
\item[$(i)$] For all $y \in Y$, $\beta>0$, $\Argmax_{x \in X} \tilde{\pi}(x, y,\beta)$ reduces to one point and defines a continuous map $\mathbf{br}$ from $Y \times \R_+^*$ to $Int(X)$.
\end{itemize}
The map $(y,\beta) \in Y \times \mathbb{R}_+^* \mapsto \mathbf{br}(y,\beta)$ is usually called a  \emph{smooth best response} map.

\hop In our analysis, we also need a little more regularity on the smooth best response map, namely we need the following: 
\begin{itemize} 
\item[$(ii)$] There exists $L>0$ such that the map $y \mapsto  \mathbf{br}(\beta,y)$ is Lipschitz continuous, with Lipschitz constant $L \beta$.
\end{itemize}

Hence we replace assumption $(A1)$ by a slightly stronger statement:
\begin{itemize}
\item[$(A1*)$] There exists $\lambda >0$ such that, for any $x \in Int(X)$ and any  $h \in TX$,
\[\left<D^2\rho(x)  h, h \right> \geq \lambda \|h\|^2.\]
\end{itemize}

\hop In particular, notice that $(A1*)$ implies that $D^2\rho(x)$ is invertible and that $\sup_{x \in Int(x)} \left\|(D^2\rho(x))^{-1} \right\| \leq \frac{1}{\lambda} < \infty$. Finally, under assumptions $(A1*)$ and $(A2)$, points $(i)$ and $(ii)$ are checked.

\hop In the remaining of the paper, we assume that $\rho$ is a \emph{good} perturbation function, i.e. a function  verifying properties $(A1*)$ and $(A2)$.

\begin{remark} Let $\rho:x \in X \mapsto \rho(x) = -\sum_{i \in I} x_i \log x_i$ be the entropy function. It is a particular case of a good perturbation function, and the resulting smooth best response is the so-called \emph{logit map}, given by
\[\mathbf{L}(\beta,y)_i = \frac{\exp (\beta \pi(i,y))}{\sum_{k \in I} \exp (\beta \pi(k,y))}\]
\end{remark}

\begin{definition} Player $1$ plays accordingly to a smooth fictitious play strategy, with the parameter $\beta>0$ (SFP($\beta$)) if
\[\mathbb{P} \left(i_{n+1}=i \mid \mathcal{F}_n \right) = \mathbf{br}(\overline{y}_n,\beta)_i, \; \, \forall n \geq 1.\]
\end{definition}

%\begin{Remarque} Asymptotically, the behaviour of $\mathbf{br}(\overline{y}_n,\beta)$ and $\mathbf{br}(\gamma_n,\beta)$ is similar, for any given prior $\overline{y}_0$. Hence, it is not difficult to prove that the asymptotic behaviour $e_n$ is the same whether player $1$ plays  a smooth best response to $\overline{y}_n$ or to $\gamma_n$.

%Expliquer pourquoi on regarde juste $\mathbf{br}(\overline{y}_n,\beta)$ et pas  $\mathbf{br}(\gamma_n,\beta)$, avec $\gamma_n = (1-1/(n+1))\overline{y}_n + 1/(n+1) \overline{y}_0$, pour un certain prior. La fonction $\Pi$ est Lipschitz et donc regarder le regret avec $\gamma_n$ est pareil asymptotiquement que le regarder avec $\overline{y}_n$.
%\end{Remarque}

\begin{theorem}[Fudenberg and Levine, 1995] \label{th:fudlev} For any $\eta > 0$, there exists $\mathbf{\beta}_0 > 0$ such that  a SFP($\beta$) strategy is $\eta$-consistent for any $\beta > \mathbf{\beta}_0$.
\end{theorem}

Smooth fictitious play is closely related to the so-called \emph{exponential weight algorithm}  and also to the \emph{follow the perturbed leader algorithm} (see \cite{CBL06}, chapters 4.2 and 4.3), even if the link with the latter is less obvious.
In \cite{HofSorVio09}, the authors discuss the consistency of continuous-time versions of FP and SFP.
%%%%%%%%%%%%%%%%%%%%%%%%%%%%%%%%%%%%%%%%%%%%%%%%%%%%%%%%%%%%%%%%%%%%%%
\subsection{Vanishingly smooth fictitious play}
%%%%%%%%%%%%%%%%%%%%%%%%%%%%%%%%%%%%%%%%%%%%%%%%%%%%%%%%%%%%%%%%%%%%%%%
A related natural strategy is given by the following. Recall that $\mathbf{br}$ is a smooth best response function, induced by a good perturbation function.

\begin{definition} Let  $\left(\mathbf{\beta}_n \right)_n$ be a sequence going to infinity. The \emph{vanishingly smooth fictitious play strategy} induced by $\beta_n$ (and $\mathbf{br}$)  for player $1$ is defined by
\[ \mathbb{P} \left( i_{n+1} = i\mid \mathcal{F}_n\right) = \mathbf{br} \left( \overline{y}_n,\mathbf{\beta}_n\right)_i \; \, \forall n \geq 1.\]
\end{definition}

\hop We use the notation VSFP($\mathbf{\beta}_n$) in the sequel. Consistency is not verified for any choice of $(\mathbf{\beta}_n)_n$. If this sequence increases too fast, then consistency might fail to hold, as shown by the following example.
\begin{example}
{\rm
Assume that, once again the game is $2$-player matching pennies and that nature uses the deterministic strategy described in example \ref{ex:1}. Then, if player one plays accordingly to a VSFP strategy induced by the logit map,  $\mathbf{\beta}_n = n$ and prior $\overline{y}_0 = (1/3,2/3)$, we have
\[\gamma_{2n} = \left(\frac{1}{2} - \frac{1}{6(2n + 1)}, \frac{1}{2} + \frac{1}{6(2n + 1)}\right) \; \mbox{ and } \; \gamma_{2n+1} = \left(\frac{1}{2} + \frac{1}{6(n + 1)}, \frac{1}{2} - \frac{1}{6(n + 1)} \right).\]
After a few lines of calculus (left to the reader) one gets:
\[\mathbb{E} \left(\delta_{l_{2n+1}} \mid \mathcal{F}_n \right) = \mathbf{L}(\gamma_{2n},\mathbf{\beta}_{2n}) = \left(\frac{1}{1+\exp \left(\frac{2n}{3(2n+1)}\right)}, \frac{1}{1+\exp \left(-\frac{2n}{3(2n+1)}\right)} \right). \]
Hence $(\pi(i_{2n+1},l_{2n+1}))_n$ is a sequence of independent random variables taking values in $\{0,1\}$, such that
\[\lim_n \mathbb{P} \left(\pi(i_{2n+1},l_{2n+1}) = 1 \right) = \lim_n \frac{1}{1+e^{2n/3(2n+1)}} = \frac{1}{1+e^{1/3}} = 1/2 - a\] with $a > 0.$
Similarly, $(\pi(i_{2n},l_{2n}))_n$ is a sequence of independent random variables taking values in $\{0,1\}$ and
\[\lim_n \mathbb{P} \left(\pi(i_{2n},l_{2n}) = 1 \right)  = \frac{1}{1+e^{2/3}} = 1/2-b\] with $b > 0.$}
\end{example}

Therefore, consistency is not satisfied for VSFP strategies with $\mathbf{\beta}_n =n$ since $e_n \rightarrow (a+b)/2 > 0$. 
\zdeux

\hop We now can state our main result
%%%%%%%%%%%%%%%%%%%%
\begin{theorem} \label{th:consistency} Any VSFP$(\mathbf{\beta}_n)$ strategy, with $\mathbf{\beta}_n \leq n^{\nu}$ for some $\nu <1$, is consistent.
\end{theorem}
%%%%%%%%%%%%%%%%%%
In \cite{BHS2}, the authors prove the same result as Theorem \ref{th:fudlev} using stochastic approximation methods. Specifically, they consider the state variable $(\overline{x}_n,\overline{y}_n,\overline{\pi}_n)_n$, write it as a stochastic approximation process relative to some differential inclusion, and prove that it almost surely converges to the \emph{consistency set}:
\[\left\{(x,y,\pi): \; \, \Pi(y) - \pi \leq \eta \right\}.\]
This is the approach taken in this paper. In section $2$ we show how our state variable can be written as a stochastic approximation algorithm, relative to some  non-autonomous differential inclusion.  A concept of Lyapunov function with respect to a set $A$ for non-autonomous systems is introduced in section $3$ and, in Proposition \ref{pr:limitsetdeterm}, we establish that $A$  attracts the so-called \emph{perturbed solutions}, under the right conditions. In our specific case, we also prove that there exists a Lyapunov function relative to the consistency set.
The proof of our main result, Theorem \ref{th:consistency}, is given in Section $4$. It consists in  showing that $(\overline{x}_n,\overline{y}_n,\overline{\pi}_n)_n$ is almost surely a perturbed solution with good properties and applying the results of Section $3$. In the appendix, we provide some general stability results for non-autonomous differential inclusions, namely we  estimate the deviation of so-called \emph{perturbed solutions} from the set of solutions curves.

%%%%%%%%%%%%%%%%%%%%%%%%%%%%%%%%%%%%%%%%%%%%%%%%%%%%%%%%%%%%%%%%%%%
%%%%%%%%%%%%%%%%%%%%%%%%%%%%%%%%%%%%%%%%%%%%%%%%%%%%%%%%%%%%%%%%%%%

%%%%%%%%%%%%%%%%%%%%%%%%%%%%%%%%%%%%%%%%%%%%%%%%%%%%%%%%%%%%%%%%%%%%%%%%%%%%%%%%%%%%
%%%%%%%%%%%%%%%%%%%%%%%%%%%%%%%%%%%%%%%%%%%%%%%%%%%%%%%%%%%%%%%%%%%%%%%%%%%%%%%%%%%%%%%%%%%%%%%%%%%%%%%%%%%%%%%%%%%%%%%%%%%%%%%%%%%%%%%%%%%%
\section{Stochastic approximations}  \label{s:SA}

%%%%%%%%%%%%%%%%%%%%%%%%%%%%%%%%%%%%%%%%%%%%%%%%%%%%%%%%%%%%%%%%%%%%%%%%%%%%%%%%%
\subsection{A stochastic difference inclusion} \label{ss:SDI}

\hop As it was previously mentioned, we are interested in the asymptotic behavior of the state variable  $v_n := (\overline{x}_n,\overline{y}_n,\overline{\pi}_n) \in M := X \times Y \times [-\|\pi\|_{\infty},\|\pi\|_{\infty}]$, where $\|\pi\|_{\infty} := \max_{i,l} |\pi(i,l)|$. We have
\[\overline{x}_{n+1} - \overline{x}_{n} - \frac{1}{n+1} \left( \delta_{i_{n+1}} - \mathbb{E}_{\sigma}(\delta_{i_{n+1}} \mid \mathcal{F}_n)\right) = \frac{1}{n+1} \left(- \overline{x}_{n} + \mathbf{br}(\overline{y}_n,\mathbf{\beta}_n)\right).\]
Writing the analogous recursive formulas for $\overline{y}_n$ and $\overline{\pi}_n$,  we obtain that
\[v_{n+1} - v_n - \frac{1}{n+1} U_{n+1} \in \frac{1}{n+1} F_n(v_n),\]
where
\begin{itemize}
\item[$-$] the noise sequence
\[U_{n+1} = (v_{n+1}-v_n) - \mathbb{E}(v_{n+1}- v_n \mid \mathcal{F}_n)\]
is a bounded martingale difference,
\item[$-$] the set valued map $F_n$ is given by
\begin{equation} \label{eq:F_n}
F_n(x,y,\pi):= \left\{(\mathbf{br}(y,\mathbf{\beta}_n) - x, \tau - y, \pi(\mathbf{br}(y,\mathbf{\beta}_n),\tau) - \pi, \; \, \tau \in Y \right\}.
\end{equation}
\end{itemize}

%%%%%%%%%%%%%%%%%%%%%%%%%%%%%%%%%%%%%%%%%%%%%%%%%%%%%%%%%%%%%%%%%%%%%%%%%%%%%%%%%%
\subsection{Stochastic approximations relative to non-autonomous differential inclusions} \label{ss:SA}

On a more general level, let  $M \subset \R^d$ and $F:\R_+ \times M \rightrightarrows \R^d$ be a set-valued map taking values in the set of non-empty, compact, convex subsets of $\R^d$. We say that $F$ is \emph{regular} if :
\begin{itemize}
\item[$(R1)$] $s \mapsto F(t,w)$ is measurable, for each $w \in M$;
\item[$(R2)$] for any $t \in \R^+$, the map $w \mapsto F(t,w)$ has a closed graph, i.e.
\[\left\{(w,w') \in M \times M : \; \, w' \in F(t,w) \right\}\]
is closed;
\item[$(R3)$] The map $F$ is uniformly bounded, i.e., $\sup_{t,w} \sup_{w' \in F(t,w)} \|w'\| \leq \|F\|_{\infty} < + \infty$.
\end{itemize}

\hop Consider a discrete time stochastic process $(v_n)_{n}$ in $M$, defined  by the recursive formula
\begin{equation}\label{algoINC}
v_{n+1} - v_n - \gamma_{n+1} U_{n+1} \in \gamma_{n+1}F_n(v_n),
\end{equation}
where $F_n: M \rightrightarrows \R^d$ is a set-valued map, $(\gamma_n)_{n}$ is a positive sequence, decreasing to $0$ and $(U_n)_{n}$ a sequence of $M$-valued random variables defined on a probability space $(\Omega, \F,P)$. Set $\tau_n := \sum_{i=1}^n \gamma_i$ and $m(s) := \sup \{ j \mid \tau_j \leq s \}$. We make the following additional assumptions:
\begin{itemize}
\item[$(SA1)$] For all $c>0$,
\[\sum_ne^{-c/\gamma_n} < \infty, \]
\item[$(SA2)$] $(U_n)_n$ is uniformly bounded (by $\|U\|_{\infty}$) and
\[\mathbb{E} \left(U_{n+1} \mid \mathcal{F}_n \right) = 0,\]
\item[$(SA3)$] The map $F:\R_+ \times M \rightrightarrows M$, given by
\[F(t,w) := F_{m(t)}(w)\]
is regular.
\end{itemize}

\begin{definition} If the conditions $(SA1)$, $(SA2)$ and $(SA3)$ are met, we say that $(v_n)$ is a \emph{good stochastic approximation algorithm relative to $F$}.
\end{definition}

\hop Call $v(\cdot)$ the continuous time affine interpolated process induced by $(v_n)_n$ and $\overline{\gamma}(\cdot)$ (resp. $\overline{U}(\cdot)$) the piecewise constant deterministic processes induced by $(\gamma_n)_n$ (resp. $(U_n)_n$):
\[v(\tau_i + s) = v_i + s \frac{v_{i+1} - v_i}{\gamma_{i+1}} \mbox{ for }  \; \, s \in [0,\gamma_{i+1}], \; \;  \overline{\gamma}(\tau_i + s) := \gamma_{i+1} \mbox{ for } \; s \in [0, \gamma_{i+1}[,\]
and analogously for $\overline{U}$.

\begin{lemma} \label{lm:pertsol} For almost every $s \in \R_+$, $v(\cdot)$ is differentiable and we have
\[\dot{v}(s) - \overline{U}(s) \in F(s,v_{m(s)}).\]
\end{lemma}

\hop {\bfseries Proof.} We have
\[v(s) = v_{m(s)} + \frac{v_{m(s)+1}-v_{m(s)}}{\gamma_{m(s)+1}}(s-\tau_{m(s)})\]
Hence, if $s \notin \{\tau_n, \, n \in \N^*\}$, $v(\cdot)$ is differentiable and
\[\dot{v}(s) = \frac{v_{m(s)+1} - v_{m(s)}}{\gamma_{m(s)+1}}.\]
Consequently
\[\dot{v}(s) - \overline{U}(s) \in F_{m(s)}(v_{m(s)}) = F(s,v_{m(s)}).\]
$\; \; \blacksquare$

\hop In the sequel, we use the notation $\overline{v}(s) := v_{m(s)}$. Notice that $\overline{v}$ is a piecewise constant map on $\R_+$.
\ztrois

\hop Let us come back to the particular case of section \ref{ss:SDI}, where $v_n = (\overline{x}_n,\overline{y}_n,\overline{\pi}_n)$ and  $F_n$ is given by  (\ref{eq:F_n}).

\begin{lemma}  $(v_n)_n$ is a good stochastic approximation algorithm with step size $\gamma_n = 1/n$, relative to the map $F$ given by  $F(t,w) := F_{m(t)}(w)$.
\end{lemma}

\hop {\bfseries Proof.} We only need to prove that $F$ is a regular set-valued map. The fact that $F$ has non-empty compact convex values is straightforward, as well as measurability. Also, the map $F$ takes values in a compact set. Thus $F$ is uniformly bounded. Given $s \in \R_+$, we now need to check upper semi-continuity of $v \mapsto F(s,v)$, which is equivalent to $\{(v,w), \; w \in F(s,v)\}$ being closed. Let $(x_n,y_n,\pi_n)$ converge to $(x,y,\pi)$. We then have $\lim_n \mathbf{br}(y_n,\mathbf{\beta}_{m(s)}) = \mathbf{br}(y,\mathbf{\beta}_{m(s)})$. Hence,
\[\lim_n \left(\mathbf{br}(y_n,\mathbf{\beta}_{m(s)}), \tau_n, \pi\left(\mathbf{br}(y_n,\mathbf{\beta}_{m(s)}),\tau_n \right)\right) =  \left(\mathbf{br}(y,\mathbf{\beta}_{m(s)}), \tau, \pi\left(\mathbf{br}(y,\mathbf{\beta}_{m(s)}),\tau \right) \right) \in F(s,x,y,\pi).\]
$\; \; \blacksquare$

In the particular case where $F$ is an autonomous set-valued map (i.e. $F$ does not depend on $t \in \R_+$), stochastic approximation algorithms described above have been studied in \cite{BHS1} and they proved that there is a deep relationship between the asymptotic behavior of $(v_n)$ and the solutions of the autonomous differential inclusion
\[\dot{\mathbf{w}} \in F(\mathbf{w}).\]
In particular, they show that, if there exists a global attractor $A$ for the deterministic dynamics, then the limit set of $(v_n)_n$ is contained in $A$.

Unfortunately, in our case, the mean deterministic system associated to our random process $(v_n)_n$ is a non-autonomous differential inclusion, as we will see later on.

%%%%%%%%%%%%%%%%%%%%%%%%%%%%%%%%%%%%%%%%%%%%%%%%%%%%%%%%%%%%%%%%%%%
%%%%%%%%%%%%%%%%%%%%%%%%%%%%%%%%%%%%%%%%%%%%%%%%%%%%%%%%%%%%%%%%%%%
\section{Lyapunov functions relative to nonautonomous differential inclusions}

%%%%%%%%%%%%%%%%%%%%%%%%%%%%%%%%%%%%%%%%%%%%%%%%%%%%%%%%%%%%%%%%%%%%%%
\subsection{Perturbed solutions and uniform Lyapunov functions}

Let us consider the \emph{non-autonomous differential inclusion}
\begin{equation} \label{eq:NADI}
\dot{\mathbf{w}}(s) \in F(s,\mathbf{w}(s)), \; \, s \in [a,b]
\end{equation}
A map $\mathbf{w}: [a,b] \rightarrow M$ is a solution of (\ref{eq:NADI}) if it is absolutely continuous and, for almost every $s \in [0,T]$, $\; \dot{\mathbf{w}}(s) \in F(s,\mathbf{w}(s))$. The existence of solutions from any initial condition is guaranteed under various sets of assumptions, in particular for regular $F$ (see Section \ref{ss:stability} for more details)

\begin{definition} A map $v: \R_+ \rightarrow M$ is a \emph{perturbed solution} of the non-autonomous differential inclusion $\dot{\mathbf{w}}(s) \in F(s,\mathbf{w}(s))$ if there is a locally integrable function $\overline{U}: \R_+ \rightarrow \R^d$ such that
\begitem
\item[$(PS1)$] $v$ is absolutely continuous,
\item[$(PS2)$] we have
\[\Delta(t,t+T) := \sup_{h \in [0,T]} \int_t^{t+h} \overline{U}(s) ds \rightarrow_{t \rightarrow + \infty} 0,\]
\item[$(PS3)$] $\dot{v}(s) - \overline{U}(s) \in F(s,\overline{v}(s))$ for some measurable map $\overline{v}: \R_+ \rightarrow M$ such that
\[\|v(s) - \overline{v}(s)\| \leq \delta(s),\]
with $\delta(s) \downarrow_s 0$.
\end{itemize}
\end{definition}

\begin{remark} Notice that, in the autonomous case, this is Definition $(II)$ in \cite{BHS1}
\end{remark}

\begin{proposition} \label{pr:pertsol} Let $v(\cdot)$ be the continuous time affine interpolated process associated to a good stochastic approximation. Then $v$ is almost surely a perturbed solution, with $\overline{v}(s) = v_{m(s)}$ and $\delta(s) = c \overline{\gamma}(s)$ (where $c$ is some positive constant).
\end{proposition}

\hop {\bfseries Proof.} This is a direct consequence of Lemma \ref{lm:pertsol} and Proposition 4.4 in \cite{Ben99}. We will provide more details in the particular case we are interested in, in Section \ref{s:proof}.
$\; \; \blacksquare$

We now define a concept of Lyapunov function adapted to non-autonomous differential inclusions.

\begin{definition} \label{def:lyapunov} Let $A$ be a compact set in $M$ and $V$ be an open neighbourhood of $A$. A smooth map $\Phi: \R_+ \times V \rightarrow \R_+$ is called a \emph{uniform Lyapunov function} for the non-autonomous differential inclusion (\ref{eq:NADI}) with respect to $A$ if the following hold:
\begitem
\item[$a)$]  we have
\[A = \left\{w \in V: \; \exists s_n \uparrow + \infty, \; \, \lim_{n \rightarrow + \infty} \Phi(s_n,w) = 0\right\},\]
\item[$b)$] There exists two maps $\lambda: \R_+^* \rightarrow ]0,1[$ and $\varepsilon: \R_+ \times \R_+ \rightarrow \R_+$ with the property that
\[\lim_{T \rightarrow + \infty} \lambda(T) = 0 \;  \mbox{ and  }  \; \, \lim_{t \rightarrow + \infty} \varepsilon(t,T) = 0, \,  \; \forall T>0;\]
and, for any $t>0, T>0$ and any solution $\mathbf{w}$ on $[t, t +T]$, we have
\[\Phi(t+s,\mathbf{w}(t+s)) \leq \lambda(s) \Phi(t,\mathbf{w}(t)) + \varepsilon(t,T), \; \forall s \in [0,T].\]
\end{itemize}
\end{definition}
If $V=M$ then $\Phi$ is called a \emph{global uniform Lyapunov function}.

\begin{remark} \label{Rm:Phiuniform} Assumption $a)$ is checked in particular if the somewhat more explicit condition is verified:
\begin{itemize}
\item[$a')$]  there exists a continuous map $g:V \rightarrow \R+$ such that
\[A = \{w \in M: g(w)=0\}, \; \, \|g(w)-\Phi(s,w)\| \rightarrow_{s \rightarrow + \infty} 0,\]
uniformly in $w \in V$.
\end{itemize}
\end{remark}
%%%%%%%%%%%%%%%%%%%%%%%%%%%%%%%%%%%%%%%%%%%%%%%%%%
%%%%%%%%%%%%%%%%%%%%

The following lemma will be useful to prove the main result of this section, namely Proposition \ref{pr:limitsetdeterm}.

\begin{lemma} \label{lm:Phi} Let $(\Phi_k)_{k \geq k_0}$, $(\lambda_k)_{k \geq k_0}$ and $(\eta_k)_{k \geq k_0}$ be positive sequences of real numbers such that $0<\lambda_k < 1$ and
\begitem
\item[$(i)$] for any $k \geq k_0$
\[\Phi_{k+1} \leq \lambda_k \Phi_k + \eta_{k+1};\]
\item[$(ii)$] for $k \geq k_0+1$, denoting $H_k := \Pi_{i=k_0}^{k-1} \lambda_i$ and $\tilde{H}_k = H_k \sum_{i=0}^{k-1} H^{-1}_i \eta_i,$ we have  $\; \lim_{k \rar \infty} H_k = \lim_{k \rar \infty} \tilde{H}_k = 0.$
\end{itemize}
Then $\lim_{k \rar \infty} \Phi_k = 0$.
\end{lemma}

\hop {\bfseries Proof.} Without loss of generality, we assume that $k_0 = 0$.  A simple recursive argument yields
\[\Phi_{k} \leq H_k  \left(\Phi_{0} + \sum_{i=1}^{k} H_i^{-1} \eta_i \right)\]
and the proof is complete.
$\; \; \blacksquare$

%\begin{Remarque}

We say that $\Phi$ is \emph{uniformly Lipschitz} if there exists $L_{\Phi} >0$ such that, for any $s \geq 0$ and $w,w' \in M$,
\[\left|\Phi(s,w) - \Phi(s,w') \right| \leq L_{\Phi}\|w-w'\|.\]

We now need to define Lipschitz continuity for non-autonomous set-valued maps: call $d_H$ the Hausdorff distance, given by
\[d_H(E,F) = \max \left\{\sup_{x \in E} d(x,F), \, \sup_{y \in F} d(y,E) \right\}.\]
Recall that $d_H$ is a pseudo-metric on the set of non-empty subsets of $M$ and a metric if we restrict to the non-empty compact sets of $M$. We say that $F$ is \emph{Hausdorff continuous} if it is continuous with respect of the Hausdorff metric:
\[\lim_{t \rightarrow t_0, w \rightarrow w_0} d_H(F(t,w),F(t_0,w_0)) = 0.\]
If $F$ is Hausdorff continuous, we call it $L$-Lipschitz, for an integrable function $L:[a,b] \rightarrow \R_+$ if
\[d_H(F(t,w),F(t,w')) \leq L(t) \|w-w'\|, \; \, \mbox{for a.e.} \;  t \in [a,b], \; \forall \; w,w'\]

We now state the main result of this section. Corollary \ref{co:v} plays an important role here, as it gives upper bound for the deviation of perturbed solutions from actual solutions of the deterministic system. For convenience of the reader, we chose to postpone this technical result to Section \ref{ss:stability}.

\begin{proposition} \label{pr:limitsetdeterm} Assume that $v$ is a perturbed solution relative to a regular Lipschitz map $F$ (with $L:\R_+ \rightarrow \R_+$) and that $\Phi$ is a global uniform Lyapunov function with respect to a compact set $A$ and the differential inclusion (\ref{eq:NADI}). Assume also that there  exists a sequence of positive real numbers  $(T_k)_k$ such that
\begin{itemize}
\item[$(i)$]  $S_k := \sum_{i=1}^k T_i \rightarrow + \infty$,
\item[$(ii)$] there exists $k_0 \in \N$ and a sequence $(r_k)_k$ such that, for any $k \geq k_0$
\[R(S_k,S_{k+1}) \leq  r_k,\]
with $R$ defined by (\ref{R}) in Corollary \ref{co:v},
\item[$(iii)$] $\Phi$ is uniformly Lipschitz, with constant $L_{\Phi}$,
\item[$(iv)$] denoting $H_k := \Pi_{i=k_0+1}^{k} \lambda(T_i)$ and $\eta_k := \varepsilon(S_{k-1},T_{k}) + L_{\Phi} r_{k-1}$, we have
\[\lim_{k \rightarrow + \infty} H_k \sum_{i=1}^k H_i^{-1} \eta_i = 0.\]

\end{itemize}
Then the limit set of $v$, $\mathcal{L}((v(s))_{s >0}) := \left\{v^*: \; \exists s_n \uparrow + \infty, \, \lim_n v(s_n) = v^*  \right\}$ is contained in $A$.
\end{proposition}

\hop {\bfseries Proof.} First, by Corollary \ref{co:v}, for any $k \in \N$, there exists a solution $\mathbf{w^k}$ on $[S_{k},S_{k+1}]$ such that $\mathbf{w^k}(S_{k}) = v(S_{k})$ and
\[\sup_{s \in [S_{k},S_{k+1}]} \|v(s) - \mathbf{w^k}(s)\| \leq R(S_k,S_{k+1}).
\]
By $(ii)$ the sequence of solutions curves $(\mathbf{w^k})_{k \geq k_0}$ is such that
\[\sup_{s \in [S_{k},S_{k+1}]} \|v(s) - \mathbf{w^{k}}(s)\| \leq r_k.\]

On the other hand, by definition of $\Phi$ and $\mathbf{w^k}$,  we have
\[\Phi(S_{k+1},\mathbf{w^k}(S_{k+1})) \leq \lambda(T_{k+1}) \Phi(S_{k},\mathbf{w^k}(S_{k})) + \varepsilon(S_{k},T_{k+1}).\]
Hence, by $(iii)$ and $(iv)$,   for any $k \geq k_0$,
\begin{eqnarray*}
\Phi(S_{k+1},v(S_{k+1})) &\leq& \Phi(S_{k+1},\mathbf{w^k}(S_{k+1})) + L_{\Phi} \left\| v(S_{k+1}) - \mathbf{w^k}(S_{k+1})\right\|\\
&\leq& \lambda(T_{k+1}) \Phi(S_{k},v(S_{k})) + L_{\Phi} r_k + \varepsilon(S_{k},T_{k+1})\\
&=& \lambda(T_{k+1}) \Phi(S_{k},v(S_{k})) + \eta_{k+1}
\end{eqnarray*}
Clearly, $H_k \rightarrow 0$, by definition on $\lambda$. Calling $\Phi_k := \Phi(S_{k},v(S_{k}))$ and $\lambda_k := \lambda(T_{k+1})$ we have $\Phi_k \rightarrow 0$ by Lemma \ref{lm:Phi}. Now let $v_*$ be a limit point of $v(s)$: $v_* = \lim_n v(s_n)$, for some sequence $s_n \uparrow_n + \infty$. Call $k(n) := \sup \{k \in \N: \; S_k \leq s_n\}$. For $n$ large enough, $k(n) \geq k_0$ and
\[\Phi(s_n,v(s_n)) \leq \lambda(s_n-S_{k(n)}) \Phi(S_{k(n)},v(S_{k(n)}) + L_{\Phi}r_{k(n)} + \varepsilon(S_{k(n)},s_n-S_{k(n)}) \rightarrow_{n \rightarrow + \infty} 0.\]
We therefore have
\[\Phi(s_n,v_*) \leq \Phi(s_n,v(s_n)) + L_{\Phi}\|v_* - v_n\| \rightarrow_{n \rightarrow + \infty} 0.\]
Consequently $v_* \in A$ and the proof is complete.
$\; \; \blacksquare$

%%%%%%%%%%%%%%%%%%%%%%%%%%%%%%%%%%%%%%%%%%%%%%%%%%%%%%%%%%%%%%%%%%%%%%%%%%%%%%%%%%%
\subsection{A Lyapunov function for the differential inclusion induced by (\ref{eq:F_n})}

We now focus on the particular case of Section \ref{ss:SDI} and prove that there exists a global Lyapunov function with respect to the so-called consistency set.

\begin{theorem} \label{th:att1} Let $A = \left\{ (x,y,\pi) \in M  \mid \Pi(y) - \pi \leq 0 \right\}$. There exists a global uniform Lyapunov  function $\Phi$ relative to the compact set $A$ and the non-autonomous differential inclusion
\begin{equation} \label{eq:NADI2}
\dot{\mathbf{w}}(s) \in F(s,\mathbf{w}(s)).
\end{equation}
\end{theorem}

\hop {\bfseries Proof.} We prove that properties $a')$ and $b)$ (of respectively Remark \ref{Rm:Phiuniform} and Definition \ref{def:lyapunov})  hold. 

\hop Let $\Phi: \R_+ \times M \rightarrow \R_+$ be defined by
\[\Phi(s,x,y,\pi) = \left\{
            \begin{array}{ll}
               \tilde{\Pi}(y,\mathbf{\beta}_{m(s)}) - \pi   &\mbox{ if }  \tilde{\Pi}(y,\mathbf{\beta}_{m(s)}) \geq \pi \\
                0 &\mbox{ if } \tilde{\Pi}(y,\mathbf{\beta}_{m(s)}) < \pi.
            \end{array}
            \right.
\]
where 
\[\tilde{\Pi}: Y \times \R_+^* \rightarrow \R, \; (y,\beta) \mapsto \max_{x \in X} \tilde{\pi}(x,y,\beta) = \tilde{\pi} \left(\mathbf{br}(y,\beta),y,\beta \right).\]

\hop Notice that
\[A = \left\{(x,y,\pi): g(x,y,\pi) = 0 \right\} \; \mbox{ and } \; \,
\|g(x,y,\pi)-\Phi(s,x,y,\pi)\| \rightarrow_{s \rightarrow + \infty} 0\]
uniformly, where $g(x,y,\pi) := \max \{0, \Pi(y) - \pi\}$. Let $t$ and $T$ be positive real numbers and $\mathbf{w}(s) := (x(s),y(s),\pi(s))$ be a solution of the non-autonomous differential inclusion (\ref{eq:NADI2}) on $[t,t+T]$, such that $\pi(s) \leq \tilde{\Pi}(y(s),\mathbf{\beta}_{m(s)})$. Thus
\[\dot{y}(s) = \tau(s) - y(s), \; \, \dot{\pi}(s) = \pi (\mathbf{br}(y(s),\mathbf{\beta}_{m(s)}),\tau(s)) - \pi(s),\]
where $\tau(s) \in Y, \; \forall s$. Let
\[\Psi(s) := \Phi(s,x(s),y(s),\pi(s)) = \tilde{\pi}\left(\mathbf{br}(y(s),\mathbf{\beta}_{m(s)}),y(s),\mathbf{\beta}_{m(s)}\right) - \pi(s).\]
Recall that $\mathbf{\beta}_{m(s)}$ is piecewise constant on $[t, t+T]$. Hence, for almost every $s \in [t,t+T]$,  we have
\begin{eqnarray*}
\dot{\Psi}(s) &=&  \tilde{\pi} \left( \mathbf{br}(y(s),\mathbf{\beta}_{m(s)}), \dot{y}(s), \mathbf{\beta}_{m(s)}\right)  - \dot{\pi}(s) \\
&=& \tilde{\pi} \left( \mathbf{br}(y(s),\mathbf{\beta}_{m(s)}), \tau(s),\mathbf{\beta}_{m(s)}\right) - \tilde{\pi} \left( \mathbf{br}(y(s),\mathbf{\beta}_{m(s)}), y(s),\mathbf{\beta}_{m(s)}\right)\\
&&- \pi (\mathbf{br}(y(s),\mathbf{\beta}_{m(s)}),\tau(s)) + \pi(s) \\
&\leq& - \Psi(s) + \frac{1}{\mathbf{\beta}_{m(s)}} \rho(\mathbf{br}(y(s),\mathbf{\beta}_{m(s)})) \leq - \Psi(s) + \frac{1}{\mathbf{\beta}_{m(s)}},
\end{eqnarray*}
where we recall that $\rho$ denotes the perturbation function. The first equality is obtained using the  enveloppe theorem and the fact that $ \tilde{\pi}$ is linear in its second argument. Thus, by an application of Gronwall's lemma, we obtain
\[\Psi(t+T) \leq e^{-T}\Psi(t) + \frac{1}{\mathbf{\beta}_{m(t)}}\]

\hop Consequently, $\Phi$ is a global uniform Lyapunov function with respect to $A$, which proves the result.
$\; \; \blacksquare$

%%%%%%%%%%%%%%%%%%%%%%%%%%%%%%%%%%%%%%%%%%%%%%%%%%%%%%%%%%%%%%%%%%%%%%%%%%%%%%%%%%%
%%%%%%%%%%%%%%%%%%%%%%%%%%%%%%%%%%%%%%%%%%%%%%%%%%%%%%%%%%%%%%%%%%%%%%%%%%%%%%%%%%%%%%%%%%%%%%%%%%%%%%%%%%%%%%%%%%%%%%%%%%%%%%%%%%%%%%%%%%%%
\section{Proof of Theorem \ref{th:consistency}} \label{s:proof}

We are now ready to prove our main result. We already proved that the interpolated random process induced by $(v_n)_n$ is almost surely a perturbed solution of the differential inclusion (\ref{eq:NADI2}) with $\delta(s) = c \overline{\gamma}(s)$, and that there exists a global uniform Lyapunov function with respect to
\[A = \left\{ (x,y,\pi) \in M  \mid \Pi(y) - \pi \leq 0 \right\},\]
see respectively Proposition \ref{pr:pertsol} and Theorem \ref{th:att1}. Therefore we now check that the assumptions of Proposition \ref{pr:limitsetdeterm} hold. Be aware that we have not used the particular form of the parameter sequence $(\mathbf{\beta}_n)_n$ so far. Recall that $\mathbf{\beta}_n = n^{\nu}$, for some $\nu  \in (0,1)$. 

Notice that $\gamma_n = 1/n$. Therefore we have $\tau_n \sim \log n$ and $m(s) = \mathcal{O}(e^s)$\footnote{more precisely, $\frac{e-1}{e} e^s \leq m(s) \leq e^s -1$}.
Recall that, given positive real numbers  $t$ and $T$,  $\Delta(t,t+T)$ denotes the random variable
\[\sup_{h \in [0,T]} \int_{t}^{t+h} \overline{U}(s) ds.\]
 Although the quantity $\mathbb{P} \left(\Delta(t,t+T) \geq \alpha \right)$ always vanishes under assumptions (SA1) and (SA2), we need to know a little more. The next lemma (proved in \cite{Duf97} or \cite{Ben99} for instance) gives an upper bound of this quantity.
\begin{lemma} \label{lm:Delta} There exists positive constants $C$ and $C'$ (depending on $\|U\|_{\infty}$) such that, for any $\alpha >0$,
\[\mathbb{P} \left(\Delta(t,t+T) \geq \alpha \right) \leq C \exp \left(\frac{- \alpha^2e^t}{C' T} \right).\]
\end{lemma}

%Notice that, by Lemmas \ref{lm:pertsol} and  \ref{lm:Delta} and a Borel-Cantelli argument, $v$ is almost surely a perturbed solution, with $\delta(s) = A\overline{\gamma}(s) \leq 2Ae^{-s}$, where $A = \|F\|_{\infty} + \|U\|_{\infty}$.

The set-valued map $F$ is regular and $L(\cdot)$-Lipschitz, with the same Lipschitz constant as the map $(s,y) \mapsto  \mathbf{br}(y,\mathbf{\beta}_{m(s)})$. Hence $L(s) = L\mathbf{\beta}_{m(s)}$, for some constant $L$ (see Section \ref{ss:perturb}). Hence, we can assume without loss of generality, that $L(s) \leq e^{\nu s}$ (up to choosing $\nu' > \nu$). In the next proposition, we see that assumptions $(i)$ and $(ii)$ of Proposition \ref{pr:limitsetdeterm} hold, if we choose $T_k = (\nu k)^{-1}$.

\begin{proposition} \label{pr:beta} If we choose $T_k := (\nu k)^{-1}$  there exist some constant $r >1$ with the property that, with probability one, there exists $k_0 \in \N$ such that points $(i)$ and $(ii)$ of Proposition \ref{pr:limitsetdeterm} are verified for $v$, with $r_k = k^{-r}$
\end{proposition}

\hop {\bfseries Proof.} Point $(i)$  clearly holds. We now need to prove $(ii)$. In this particular case, the quantity $R(S_k,S_{k+1})$ satisfies 
\[R(S_k,S_{k+1}) \leq \left(\Delta(S_{k},S_{k+1}) + c \overline{\gamma}(S_{k})\right) \exp \left( \int_{S_{k}}^{S_{k+1}} L(\tau) d\tau \right).\]
By our choice of the sequence $(T_k)_k$, $\exp (\nu S_k) \leq \exp \left(1+ \log k \right) \leq 3k$. Hence
\[\exp \left( \int_{S_k}^{S_{k+1}} L(\tau) d\tau \right) \leq \exp(T_{k+1} e^{\nu S_{k+1}}) \leq C_0,\]
for some constant $C_0$ which depends on $\nu$. Additionally, $\overline{\gamma}(S_k) \leq 2e^{-S_k} \leq 2 k^{-1/\nu}$. Hence
\[c\overline{\gamma}(S_k) \exp \left(T_{k+1} e^{\nu S_{k+1}}\right) \leq \frac{3c}{k^{1/\nu}}. \]
Choose $r \in (1,\frac{\nu +1}{2 \nu})$. By Lemma \ref{lm:Delta},
\begin{eqnarray*}
\mathbb{P} \left( \Delta(S_k,S_{k+1}) \exp \left( \int_{S_k}^{S_{k+1}} L(\tau) d\tau \right) \geq \frac{1}{2k^{r}} \right) &\leq& C \exp \left(\frac{- k^{-2r} e^{S_k}}{4C'C_0 T_{k+1}} \right)\\
&\leq& C \exp \left(\frac{- k^{-2r + 1/\nu}}{C'C_0 \nu^{-1} (k+1)^{-1}}\right)\\
&\leq& C \exp\left(\frac{- k^{-2r + 1 + 1/\nu}}{C'_1}\right)
\end{eqnarray*}
for some positive constant $C'_1$. Now, since $r < 1/\nu$,  we have  for $k$ large enough
\[c\overline{\gamma}(S_k) \exp \left(T_{k+1} e^{\nu S_{k+1}}\right) \leq \frac{1}{2k^{r}}. \]

\hop Consequently, if we call $A_k$ the event
\[\left\{ \left(\Delta(S_{k},S_{k+1}) + c \overline{\gamma}(S_k)\right) \exp \left( \int_{S_{k}}^{S_{k+1}} L(\tau) d\tau \right)  \geq  \frac{1}{k^{r}}\right\},\]
then
\[\mathbb{P} \left( A_k \right) \leq C \exp \left(\frac{- k^{-2r + 1 + 1/\nu}}{C'_1}\right).\]
By an application of the Borel-Cantelli lemma, with probability one, there exists $k_0 \in \N$ such that, for any $k \geq \N$,
\[R(S_k,S_{k+1}) \leq \left(\Delta(S_{k},S_{k+1}) + c \overline{\gamma}(S_{k})\right) \exp \left( \int_{S_{k}}^{S_{k+1}} L(\tau) d\tau \right)  \leq  \frac{1}{k^{r}},\]
which yields the result.
$\; \; \blacksquare$

\begin{remark} \label{rm:Ls} By similar arguments, we can also prove the following: Assume that $F$ is $L$-Lipschitz, with  $L(s) \leq L s$. Then there exist $T >0$, and $r>0$ such that, with probability one, there exists $k_0 \in \N$ with the property that points $(i)$ and $(ii)$ of Proposition \ref{pr:limitsetdeterm} are verified for $v$, with $T_k = T$ and $r_k = e^{-r k}$.
%\[\mathbb{P} \left(\inf_{\mathbf{x} \in \mathcal{S}(t,t+T)}\sup_{s \in [t,t+T]} \|v(s) - \mathbf{x}(s)\| \geq \varepsilon \right) \leq C \exp \left(\frac{- \varepsilon^2 \exp \left(\gamma t \right)}{4C'T}\right),\]
%for some constant $\gamma>0$.
\end{remark}

 Consequently,  points $(i)$ and $(ii)$ of Proposition \ref{pr:limitsetdeterm} are  almost surely satisfied for $k \geq k_0$, with $T_k = (\nu k)^{-1}$ and $r_k = k^{-r}, \; r>1$. We now need to check points $(iii)$ and $(iv)$.

\hop Let $b$ be a positive constant and consider the map $\phi:Y \times [-\|\pi\|_{\infty},\|\pi\|_{\infty}] \rightarrow \R_+$, given by
\[\phi(y,\pi) = \left\{
            \begin{array}{ll}
               \tilde{\Pi}(y,b) - \pi   &\mbox{ if }  \tilde{\Pi}(y,b) \geq \pi \\
                0 &\mbox{ if } \tilde{\Pi}(y,b) < \pi.
            \end{array}
            \right.
\]
Let $(y,\pi)$ be such that $\tilde{\Pi}(y,b) > \pi$. Then, by Lemma 6.2 in \cite{BHS2} (see also \cite{FudLev99}), we have
\begin{eqnarray*}
\frac{\partial}{\partial y} \phi(y,\pi)(h) &=& \frac{\partial}{\partial y} \tilde{\pi} (\mathbf{br}(y,b),y,b) (h)\\
&=& \pi(\mathbf{br}(y,b),h)
\end{eqnarray*}
and
\[\frac{\partial}{\partial \pi} \phi(y,\pi) = -1.\]
Thus
\[|\phi(y,\pi) - \phi(y',\pi')| \leq \|\pi\|_{\infty} \|y-y'\| + |\pi - \pi'|\]
and $\phi$ is Lipschitz with Lipschitz constant that does not depend on $b$, which means that the map $v \mapsto \Phi(s,v)$ is uniformly Lipschitz.

We now prove point $(iv)$. By Theorem \ref{th:att1}, $\Phi$ is a global uniform Lyapunov function relative to
\[A= \left\{ (x,y,\pi) \in M  \mid \Pi(y) - \pi \leq 0 \right\},\]
with $\lambda(T) = e^{-T}$ and $\varepsilon(t,T) = \frac{T}{\mathbf{\beta}_{m(t)}}$. Hence \[\eta_{k+1} = L_{\Phi} k^{-r} + \frac{T_{k+1}}{\mathbf{\beta}_{m(S_k)}} \leq L_{\Phi} k^{-r}+ c \frac{T_{k+1}}{k},\]
for some positive constant $c$. We have $\sum_i \eta_i < \infty$ and $H_k = e^{-\sum_{i=k_0}^k T_i} = \mathcal{O}(k^{-1/\nu})$. Thus point (iv) is checked (see point b) of Lemma \ref{lm:Hk} for a proof).

As a consequence, Proposition \ref{pr:limitsetdeterm} applies and
\[\mathcal{L}((v(s))_{s>0}) \subset A\]
almost surely. In particular
\[\limsup_n e_n \leq 0, \; \; \mbox{almost surely}\]
and Theorem \ref{th:consistency} is proved.

%%%%%%%%%%%%%%%%%%%%%%%%%%%%%%%%%%%%%%%%%%%%%%%%%%%%%%%%%%%%%%%%%%%%%%%%%%%%%%%%%%%%%%%%%%%%%%%%%%%%%%%%%%%%%%%%%%%%%%%%%%%%%%%%%%%%%%%%%%%%
\section{Appendix}

\subsection{Sufficient conditions for Lemma \ref{lm:Phi}, $(ii)$ to hold}

\begin{lemma}
 \label{lm:Hk} Point $(ii)$ of Lemma \ref{lm:Phi} is verified in the following cases:
\begitem
\item[$a)$] $\lambda_k = \lambda <1$ and $\lim_{k \rar \infty} \eta_k = 0$,
\item[$b)$] $\lim_{k \rar \infty} H_k =  0$ and $\sum_i \eta_i < +\infty$.
\end{itemize}
\end{lemma}
%\end{Remarque}

\hop {\bfseries Proof.} For point $a)$, $H_k = \lambda^k$ and we have
\begin{eqnarray*}
\tilde{H}_{k+k'} &=& \lambda^{k+k'} \left(\sum_{i=1}^k H_i^{-1} \eta_i + \sum_{i=k+1}^{k+k'} H_i^{-1}\eta_i \right)\\
&\leq&  \lambda^{k'} \max_{i=1,...,k} \eta_i + \eta_{k+1} \sum_{i=0}^{k'-1} \lambda^{i}\\
&\leq& \lambda^{k'} \max_{i=0,...,k} \eta_i + \eta_{k+1} \frac{1}{1-\lambda},
\end{eqnarray*}
which gives the result.

\hop For the second point, remember that $(H_k)_k$ is a decreasing sequence. Hence
\begin{eqnarray*}
\tilde{H}_{k+k'} &\leq& H_{k+k'} \left(\sum_{i=1}^k H_i^{-1} \eta_i + H_{k+k'}^{-1} \sum_{i=k+1}^{k+k'} \eta_i \right)\\
&\leq& H_{k+k'} \left(\sum_{i=1}^{k} H^{-1}_i \eta_i\right) + \sum_{i=k+1}^{+ \infty} \eta_{i}.
\end{eqnarray*}
Given $\varepsilon >0$, by choosing $k$ large enough, the second term is smaller than $\varepsilon$. Then we can pick $k'$ large enough so that the first term is also smaller than $\varepsilon$ and the proof is complete.
$\; \; \blacksquare$

\subsection{Stability of one-sided Lipschitz differential inclusions} \label{ss:stability}

Let  $M \subset \R^d$. Consider a set-valued map $F:\R_+ \times M \rightrightarrows M$ taking values in the set of non-empty, compact, convex subsets of $M$. Given $a<b$, let us consider the non-autonomous differential inclusion (\ref{eq:NADI}):
\begin{equation} 
\dot{\mathbf{w}}(s) \in F(s,\mathbf{w}(s)), \; \, s \in [a,b].
\end{equation}
For $A \subset M$ we let  $F^{-1}(A) = \{(s,w) \in [a,b] \times M: \; \, F(s,w) \cap A \neq \emptyset\}.$
We say that $F$ is measurable if $F^{-1}(A)$ is measurable, for any closed set $A \subset M$. It is \emph{upper semi-continuous} (USC) (resp. \emph{lower semi-continuous} (LSC)) if, for any closed (resp. open) set $A \subset M$, $F^{-1}(A)$ is closed (resp. open) in $[a,b] \times M$. If $M$ is compact, $F$ is upper semi-continuous if and only if its graph is closed.

We now introduce a regularity condition:
\begin{definition}[Relaxed One-sided Lipschitz] we say that the set-valued map $F$ is \emph{Relaxed One-sided Lipschitz} (ROSL) on $[a,b] \times M$ if there exists an integrable map $L:[a,b] \rightarrow M$ such that, for any $t,t'$ in $[a,b]$ $w,w' \in M$ and any $y \in F(t,w)$ there exists $y' \in F(t',w')$ with
\[<w'-w \mid y' -y> \leq L(t) \|x'-x\|^2, \; \, \forall  t \in [a,b].\]
\end{definition}

\begin{remark} If $F$ is $L(\cdot)$-Lipschitz then it is $L(\cdot)$-ROSL.
\end{remark}

The question of existence of solutions to (\ref{eq:NADI}) has been studied extensively. One of the first result on the topic was proved by Filippov (see \cite{Fil72}) and says that if $F(\cdot,\cdot)$ is Hausdorff continuous on any closed set of $[a,b] \times M$ then, for any $w_0 \in M$, there exists a solution $\mathbf{w}(\cdot)$ of (\ref{eq:NADI}), with $\mathbf{w}(a) = x_0$. Under less restrictive assumptions, the same result still holds (see \cite{Ole75}; on the topic, see also \cite{HimVan86}).

\begin{theorem}[Olech, 1975] \label{th:olech} Assume that $F$ is regular. Then there exists a solution  $\mathbf{w}(\cdot)$ of (\ref{eq:NADI}), with $\mathbf{w}(a) = w_0$.
\end{theorem}

\hop The following result will prove useful to establish Theorem \ref{th:closesolution}.

\begin{lemma} \label{co:gron2} Let $y$ be a continuously differentiable function on $[a,b]$ and $f$, $g$ be non-negative, continuous maps. If, for every $s \in [a,b]$, $\|\dot{y}(s)\| \leq f(s) \|y(s)\| + g(s)$ then
\[\|y(s)\| \leq \|y(a)\| \exp \left(\int_a^s f(\tau) d \tau \right) + \int_a^s g(u) \exp \left(\int_u^s f(\tau) d \tau \right) ds\]
\end{lemma}

\hop {\bfseries Proof.} Notice that
\[\|y(s)\| \leq \|y(a)\| + \int_a^s \|\dot{y}(u)\| du \leq \|y(a)\| + \int_a^s g(u) du + \int_a^s f(u) \|y(u)\|  du\]
and apply the integral form of Gronwall's lemma.
$\; \; \blacksquare$

In the remaining of this section, we assume that $F$ is regular. The set of solution trajectories on $[a,b]$ (resp. starting in $w_0$) will be labelled $\mathcal{S}(a,b)$ (resp. $\mathcal{S}(w_0,a,b)$).

%\begin{Remarque} Existence of solutions still holds under the weaker condition
%\begin{equation} \label{eq:intbounded} \sup_{y \in F(s,x)} \|y\| \leq c(s) (1 + \|x\|),\end{equation}
%where $c(\cdot)$ is an integrable positive function on $I$. However, we need measurability of $(t,x) \mapsto F(t,x)$ in this case.
%\end{Remarque}

\begin{theorem} \label{th:closesolution} Let $W: [a,b] \rightarrow M$ be an absolutely continuous function such that there exists a measurable map $\overline{v}: [a,b] \rightarrow M$ and a bounded measurable map $r:[a,b] \rightarrow \R_+$ which satisfy, for almost every $s \in [a,b]$,
\[d(W(s),\overline{v}(s)) \leq r(s), \; \; \dot{W}(s) \in F(s,\overline{v}(s)).\]
Then
\begin{itemize}
\item[$a)$] if $F$ is ROSL with respect to the integrable function $L$, then there exists a solution $\mathbf{w}:[a,b] \rightarrow M$ of (\ref{eq:NADI}) such that $\mathbf{w}(a) = W(a)$ and
\[\sup_{s \in [a,b]} \|\mathbf{w}(s) - W(s)\|^2 \leq \int_{a}^b  \alpha(s)  \exp \left( 4\int_{s}^b L(\tau) d \tau \right)ds,\]
where $\alpha(s) = 4 L(s) r^2(s) + 4 r(s) \|F\|_{\infty}$.
\item[$b)$] if we now assume that $F$ is Lipschitz continuous, with respect to $L$ then the conclusions of $a)$ trivially still hold and $\mathbf{w}$ can also be chosen such that
\[\sup_{s \in [a,b]} \|\mathbf{w}(s) - W(s)\| \leq \int_a^b r(s) L(s) \exp \left(\int_s^b L(\tau) d \tau \right) ds.\]
\end{itemize}
\end{theorem}

\hop {\bfseries Proof.} We prove the first point. Consider the set-valued map $G:[a,b] \times M \rightrightarrows M$ given by
\[G(s,x) := \left\{v \in F(s,w): \; \, (w-W(s) \mid v-\dot{W}(s)) \leq 2L(s) \|w-W(s)\|^2 + \frac{1}{2}\alpha(s) \right\}.\]
For any $(s,w)$, the set $G(s,w)$ is non-empty. Indeed, by the ROSL condition, since $\dot{W}(s) \in F(s,\overline{v}(s))$, there exists $v \in F(s,w)$ such that
\[(w-\overline{v}(s) \mid v-\dot{W}(s)) \leq L(s) \|w-\overline{v}(s)\|^2.\]
Hence we have
\begin{eqnarray*}
(w-W(s) \mid v-\dot{W}(s)) &\leq& L(s) \|w-\overline{v}(s)\|^2 + \|\overline{v}(s)-W(s)\| ( \|v\| + \|\dot{W}(s)\|) \\
&\leq&  2L(s) \|W(s)-w\|^2 + 2L(s)r(s)^2 + 2r(s)\|F\|_{\infty}\\
&=& 2L(s) \|W(s)-w\|^2 + \frac{1}{2} \alpha(s).
\end{eqnarray*}
Now clearly, the set $G(s,w)$ is compact and convex. The map $w \mapsto G(s,w)$ has a closed graph, for any $s \in [a,b]$. Finally It is measurable in $s$ since every map involved is measurable. Consequently, there exists a solution to the non-autonomous differential inclusion
\[\dot{\mathbf{w}}(s) \in G(s,\mathbf{w}(s)),\]
with initial condition $\mathbf{w}(a)=W(a)$. In particular, $\mathbf{w}$ is a solution of (\ref{eq:NADI}) and we also have, for almost every $s$
\[(\mathbf{w}(s)-W(s) \mid \dot{\mathbf{w}}(s) -\dot{W}(s)) \leq 2L(s) \|W(s)-\mathbf{w}(s)\|^2 + \frac{1}{2} \alpha(s).\]
Hence , for almost every $s$, we have
\begin{eqnarray*}
\frac{d}{ds} \|\mathbf{w}(s)-W(s)\|^2 &=& 2 (\mathbf{w}(s)-W(s) \mid \dot{\mathbf{w}}(s) -\dot{W}(s)) \\
&\leq& 4 L(s) \|W(s)-\mathbf{w}(s)\|^2 + \alpha(s)
\end{eqnarray*}
and point $a)$ follows from the differential form of  Gronwall's lemma.

When the Lipschitz continuity holds, let us consider the set-valued map $H: [a,b] \times M \rightrightarrows M$ given by
\[H(s,w) := \left\{v \in F(s,w): \; \, \|v-\dot{W}(s)\| \leq L(s) \|w-W(s)\| + L(s) r(s) \right\}.\]
The fact that $H$ has non-empty values follows from Lipschitz continuity: given $s$ and $w$, since $\dot{W}(s) \in F(s,\overline{v}(s))$, there exists $v \in F(s,w))$ such that
\[\|v-\dot{W}(s)\| \leq L(s) \|w-\overline{v}(s)\| \leq L(s) \left( \|w-W(s)\| + \|W(s)-\overline{v}(s)\|\right).\]
Hence $v \in H(s,w) \neq \emptyset$. Also $H(s,w)$ is convex and compact, the map $w \mapsto H(s,w)$ has a closed graph and $s \mapsto H(s,w)$ is measurable. Thus, there exists a solution $\mathbf{x}$ to the non-autonomous differential inclusion
\[\dot{\mathbf{w}}(s) \in H(s,\mathbf{w}(s)),\]
with initial condition $\mathbf{w}(a)=W(a)$. In particular, $\mathbf{w}$ is a solution of (\ref{eq:NADI}) and we also have, for almost every $s$
\[\|\dot{\mathbf{w}}(s)-\dot{W}(s)\| \leq L(s) \|\mathbf{w}(s)-W(s)\| + L(s) r(s)\]
By Gronwall's lemma (see Lemma \ref{co:gron2}), we then have
\[\sup_{s \in [a,b]} \|\mathbf{w}(s)-W(s)\| \leq \int_a^b L(s)r(s) \exp \left(\int_s^b L(\tau) d \tau \right) ds\]
and point $b)$ is proved.
$\; \; \blacksquare$

\begin{corollary} \label{co:v} Let $v:[a,b] \rightarrow M$ be an absolutely continuous map. Assume that there exist  measurable maps $\overline{v}: [a,b] \rightarrow M$,  $\delta:[a,b] \rightarrow \R_+$ bounded  and $\overline{U}:[a,b] \rightarrow M$ integrable such that, for almost every $s \in [a,b]$,
\[\dot{v}(s) - \overline{U}(s) \in F(s,\overline{v}(s)), \; \, \|v(s) - \overline{v}(s)\| \leq \delta(s).\]
Then
if $F$ is $L(\cdot)$-Lipschitz, there exists a solution $\mathbf{w}$ on $[a,b]$ such that $\mathbf{w}(a) = v(a)$ and
\[\sup_{s \in [a,b]} \|v(s) - \mathbf{w}(s)\| \leq R(a,b),\]
where
\begin{equation} \label{R} R(a,b) = \Delta(a,b) \exp \left(\int_{a}^b L(\tau) d \tau \right) + \sup_{s \in [a,b]} \delta(s) \left(\exp \left(\int_{a}^b L(\tau) d \tau \right)-1 \right)
\end{equation}
and $\Delta(a,b) = \sup_{s \in [a,b]} \|\int_a^s \overline{U}(\tau) d \tau\|$.
\end{corollary}

\hop {\bfseries Proof.} Define $W:[a,b] \rightarrow M$ by
\[W(s) :=  v(s) - \int_a^s \overline{U}(\tau) d\tau.\]
Clearly, $W$ is absolutely continuous and, for any $s$ for which $v$ is differentiable, we have $\dot{W}(s) = \dot{v}(s) - \overline{U}(s) \in F(s,\overline{v}(s))$. Additionally,
\[\|W(s) - \overline{v}(s)\| \leq \|v(s)-\overline{v}(s)\| + \|\int_a^s \overline{U}(\tau) d\tau\| \leq \delta(s) + \left\|\int_a^s \overline{U}(\tau) d\tau \right\|.\]
By a direct application of Theorem \ref{th:closesolution} with
$r(s) = \delta(s) + \left\|\int_a^s \overline{U}(\tau) d\tau \right\|,$ there exists a solution $\mathbf{w}$ such that $\mathbf{w}(a) = W(a) = v(a)$ and
%when $F$ is ROSL($L(\cdot)$),
%\[.\]
%If $F$ is $L(\cdot)$-Lipschitz, we obtain
\begin{eqnarray*}
\sup_{s \in [a,b]} \|v(s) - \mathbf{w}(s)\| &\leq& \Delta(a,b) + \int_a^b L(s) \left(\delta(s) + \|\int_a^s \overline{U}(\tau) d \tau\|\right)  \exp \left(\int_s^b L(\tau) d\tau \right)ds \\
&\leq& \Delta(a,b) + (\sup_{s \in [a,b]} \delta(s) + \Delta(a,b)) \int_a^b L(s)\exp \left(\int_s^b L(\tau) d\tau \right) ds  \leq R(a,b).
\end{eqnarray*}
$\; \; \blacksquare$

% Acknowledgments here
\section*{Acknowledgments.}
We acknowledge financial support from the Swiss National Science Foundation Grant 200020-130574. We thank Drew Fudenberg and Satoru Takahashi for suggesting us to work on this question. We also thank David Leslie for his careful reading, as well as his useful comments.

% References here (outcomment the appropriate case) 

% CASE 1: BiBTeX used to constantly update the references 
%   (while the paper is being written).
%\bibliographystyle{ormsv080} % outcomment this and next line in Case 1
%\bibliography{<your bib file(s)>} % if more than one, comma separated

% CASE 2: BiBTeX used to generate mypaper.bbl (to be further fine tuned)
%\input{mypaper.bbl} % outcomment this line in Case 2

\bibliographystyle{ormsv080}
\bibliography{MOR}

\begin{thebibliography}{23}
\expandafter\ifx\csname natexlab\endcsname\relax\def\natexlab#1{#1}\fi
\expandafter\ifx\csname url\endcsname\relax
  \def\url#1{{\tt #1}}\fi
\expandafter\ifx\csname urlprefix\endcsname\relax\def\urlprefix{URL }\fi
\expandafter\ifx\csname urlstyle\endcsname\relax
  \expandafter\ifx\csname doi\endcsname\relax
  \def\doi#1{doi:\discretionary{}{}{}#1}\fi \else
  \expandafter\ifx\csname doi\endcsname\relax
  \def\doi{doi:\discretionary{}{}{}\begingroup \urlstyle{rm}\Url}\fi \fi

\bibitem[{Bena{\"i}m(1999)}]{Ben99}
Bena{\"i}m, M. 1999.
\newblock {Dynamics of stochastic approximation algorithms}.
\newblock {\it S{\'e}minaire de probabilit{\'e}s de Strasbourg\/} {\bf 33}
  1--68.

\bibitem[{Bena{\"i}m and Hirsch(1999)}]{BenHir99}
Bena{\"i}m, M., M.W. Hirsch. 1999.
\newblock {Mixed Equilibria and Dynamical Systems Arising from Fictitious Play
  in Perturbed Games}.
\newblock {\it Games and Economic Behavior\/} {\bf 29}(1-2) 36--72.

\bibitem[{Bena{\"i}m et~al.(2005)Bena{\"i}m, Hofbauer, and Sorin}]{BHS1}
Bena{\"i}m, M., J.~Hofbauer, S.~Sorin. 2005.
\newblock {Stochastic approximations and differential inclusions}.
\newblock {\it I. SIAM Journal on Optimization and Control\/} {\bf 44}
  328--348.

\bibitem[{Bena{\"i}m et~al.(2006)Bena{\"i}m, Hofbauer, and Sorin}]{BHS2}
Bena{\"i}m, M., J.~Hofbauer, S.~Sorin. 2006.
\newblock {Stochastic Approximations and Differential Inclusions. Part II:
  Applications}.
\newblock {\it Mathematics of Operations Research\/} {\bf 31} 673--695.

\bibitem[{Blackwell(1954)}]{Bla54}
Blackwell, D. 1954.
\newblock Controlled random walks.
\newblock {\it Proceedings of the International Congress of Mathematicians\/},
  vol.~3. 336--338.

\bibitem[{Brown(1951)}]{Bro51}
Brown, G.W. 1951.
\newblock {Iterative solution of games by fictitious play}.
\newblock {\it Activity analysis of production and allocation\/} {\bf 13}(1)
  374--376.

\bibitem[{Cesa-Bianchi and Lugosi(2006)}]{CBL06}
Cesa-Bianchi, N., G.~Lugosi. 2006.
\newblock {\it {Prediction, learning, and games}\/}.
\newblock Cambridge Univ Pr.

\bibitem[{Duflo(1997)}]{Duf97}
Duflo, M. 1997.
\newblock {\it {Random iterative models}\/}.
\newblock Springer Verlag.

\bibitem[{Filippov(1971)}]{Fil72}
Filippov, AF. 1971.
\newblock {The existence of solutions of generalized differential equations}.
\newblock {\it Mathematical Notes\/} {\bf 10}(3) 608--611.

\bibitem[{Foster and Vohra(1993)}]{FosVoh93}
Foster, D.P., R.V. Vohra. 1993.
\newblock A randomization rule for selecting forecasts.
\newblock {\it Operations Research\/}  704--709.

\bibitem[{Foster and Vohra(1998)}]{FosVoh98}
Foster, D.P., R.V. Vohra. 1998.
\newblock Asymptotic calibration.
\newblock {\it Biometrika\/} {\bf 85}(2) 379--390.

\bibitem[{Fudenberg and Kreps(1993)}]{FudKre93}
Fudenberg, D., D.~Kreps. 1993.
\newblock {Learning mixed equilibria}.
\newblock {\it Games and Economic Behavior\/} {\bf 5}(3) 320--367.

\bibitem[{Fudenberg and Levine(1995)}]{FudLev95}
Fudenberg, D., D.K. Levine. 1995.
\newblock {Consistency and cautious fictitious play}.
\newblock {\it Journal of Economic Dynamics and Control\/} {\bf 19}(5-7)
  1065--1089.

\bibitem[{Fudenberg and Levine(1998)}]{FudLev98}
Fudenberg, D., D.K. Levine. 1998.
\newblock {\it {The Theory of Learning in Games}\/}.
\newblock MIT Press.

\bibitem[{Fudenberg and Levine(1999)}]{FudLev99}
Fudenberg, D., D.K. Levine. 1999.
\newblock {Conditional Universal Consistency* 1}.
\newblock {\it Games and Economic Behavior\/} {\bf 29}(1-2) 104--130.

\bibitem[{Hannan(1957)}]{Han57}
Hannan, J. 1957.
\newblock {Approximation to Bayes risk in repeated play}.
\newblock {\it Contributions to the Theory of Games\/} {\bf 3} 97--139.

\bibitem[{Harsanyi(1973)}]{Har73}
Harsanyi, J.C. 1973.
\newblock {Games with randomly disturbed payoffs: A new rationale for
  mixed-strategy equilibrium points}.
\newblock {\it International Journal of Game Theory\/} {\bf 2}(1) 1--23.

\bibitem[{Hart and Mas-Colell(2001)}]{HarMas01}
Hart, S., A.~Mas-Colell. 2001.
\newblock A general class of adaptive strategies.
\newblock {\it Journal of Economic Theory\/} {\bf 98}(1) 26--54.

\bibitem[{Himmelberg and Van~Vleck(1986)}]{HimVan86}
Himmelberg, CJ, FS~Van~Vleck. 1986.
\newblock {Existence of solutions for generalized differential equations with
  unbounded right-hand side* 1}.
\newblock {\it Journal of Differential Equations\/} {\bf 61}(3) 295--320.

\bibitem[{Hofbauer and Sandholm(2002)}]{HofSan02}
Hofbauer, J., W.H. Sandholm. 2002.
\newblock {On the Global Convergence of Stochastic Fictitious Play}.
\newblock {\it Econometrica\/} {\bf 70}(6) 2265--2294.

\bibitem[{Hofbauer et~al.(2009)Hofbauer, Sorin, and Viossat}]{HofSorVio09}
Hofbauer, J., S.~Sorin, Y.~Viossat. 2009.
\newblock Time average replicator and best-reply dynamics.
\newblock {\it Mathematics of Operations Research\/} {\bf 34}(2) 263--269.

\bibitem[{Olech(1975)}]{Ole75}
Olech, C. 1975.
\newblock {Existence of solutions of non-convex orientor fields}.
\newblock {\it Boll. Unione Mat. Ital.\/} {\bf 11} 189--197.

\bibitem[{Perchet(2010)}]{Per10}
Perchet, V. 2010.
\newblock {Approchabilit{\'e}, Calibration et Regret dans les Jeux {\`a}
  Informations Partielles}.
\newblock Ph.D. thesis, UPMC.

\end{thebibliography}

\end{document}